\input amstex
\input epsf
\documentstyle{amsppt}
\magnification\magstep1
\NoBlackBoxes
\loadmsbm
\pageheight{7.5in}
\pagewidth{5.5in}
\def\sbs{\subset}
\def\sbsq{\subseteq}

\def\es{\emptyset}
\def\sm{\setminus}

\def\i{\text{int}}
\def\t{\text}
\def\a{\alpha}
\def\e{\epsilon}
\def\O{\Omega}
\def\s{\sigma}

\topmatter
\author
J.E. KEESLING AND C. KRISHNAMURTHI
\endauthor

\title
THE SIMILARITY BOUNDARY OF A SELF-SIMILAR SET
\endtitle

\abstract 
We define the similarity boundary of a self-similar set and use it to analyze
the properties of self-similar sets in the general setting of any complete
metric space. The similarity boundary is an attempt at extending the concept of
the topological boundary in a way that is consistent with our intuitive
understanding of the term boundary. We also show how with the analysis of the
similarity boundary, we can restrict ourselves to the space $K$ and disregard
the ambient space.
\endabstract

\address 
University of Florida, Department of Mathematics, P.O. Box 118105,
358 Little Hall, Gainesville, FL 32611-8105, USA
\endaddress

\email 
jek\@math.ufl.edu
\endemail

\address 
University of Florida, Department of Mathematics, P.O. Box 118105,
358 Little Hall, Gainesville, FL 32611-8105, USA
\endaddress

\email 
cak\@math.ufl.edu
\endemail

\subjclass 
Primary 28A78, 54E40; Secondary 54H15
\endsubjclass

\keywords  Hausdorff dimension,
iterated function systems, $\beta$-space, self-similar fractals,
image measure, sequence space, sef-similar tiles.
\endkeywords

\endtopmatter

\document

\heading
\S1. Introduction
\endheading
 
We analyze self-similar sets which arise as the invariant sets or attractors
of a finite collection of similitudes on a complete metric space. Let 
$\{f_1, \cdots, f_N\}$ be a collection of contracting similitudes with 
contraction factors $\{r_1, \cdots, r_N\}$, defined on a complete metric 
space $(X, d)$. Then there exists a unique, non-empty, compact subset 
$K \sbsq X$ such that $K = \bigcup_{i=1}^N f_i(K)$. This set $K$, which is the
attractor of the collection of similitudes is called the self-similar set or
the invariant set of the similitudes. The collection of similitudes is called
an iterated function system (IFS). 

There are a number of dimensions associated
with a self-similar set. Of these, one of the most important and widely studied
dimensions is the Hausdorff dimension and the associated Hausdorff measure. In 
fact, one of the earliest characterizations of self-similar sets with positive
Hausdorff measure gave rise to the famous {\it open set condition} (OSC) set 
down by Hutchinson[8]. If there exists a non-empty, bounded, open set $G$ in 
${\Bbb R}^n$ such that $f_i(G) \subset G$ and 
$f_i(G) \cap f_j(G) = \es $ for all $i \neq j$ 
then the set of similitudes $\{f_1, \cdots, f_N \}$ is said to satisfy the OSC.
In the setting of the Euclidean metric space, the OSC is equivalent to finite
positive Hausdorff measure in dimension $\a$ where $\a$ is the similarity 
dimension (to be defined later).
Since positive Hausdorff measure implied that the Hausdorff dimension coincided
with the similarity dimension, the OSC guaranteed that the Hausdorff dimension
was $\a$. But in the general setting of any complete metric space, the OSC 
does not
even imply that the Hausdorff dimension is $\a$. This led to a strengthening
of the OSC - called the {\it strong open set condition} (SOSC) - 
attributed to Lalley. In addition to the OSC being satisfied, the SOSC requires that the open
set $G$ be such that $G \cap K \neq \es$. The SOSC implies that the Hausdorff
dimension is $\a$ but does not imply that the Hausdorff measure is positive.
Subsequently, there have been other characterizations of self-similar sets with
positive Hausdorff measures both in the setting of Euclidean space and any
complete metric space in general ([1], [15]).   

Both the OSC and the SOSC are weaker versions of the requirement
that: $f_i(K) \cap f_j(K) = \es$ for all $i \neq j$. 
When this last condition holds $\text{dim}_H K = \alpha$ and the Hausdorff 
$\alpha$-measure of $K$ is finite and positive whatever the ambient space.
One can think of the OSC and SOSC as attempts to guarantee that
$\text{dim}_H K = \alpha$ while allowing $f_i(K) \cap f_j(K)$ to be possibly
nonempty for $i \neq j$.
This has led to
serious interest in intersections of the form $f_i(K) \cap f_j(K)$
in a quest to further weaken the OSC and SOSC. In this 
paper we study these intersections through a measure on $K$ called the image 
measure (to be defined in section 2). We define what we call the similarity 
boundary of $K$ which is an attempt at extending the concept of a topological 
boundary. This approach to self-similar sets is useful because we have been 
able to restrict ourselves to the set $K$. In other words, the ambient space 
$X$ does not play a role. This is significant because, in the context of $K$,
the open set in the SOSC is required to be open only in $K$. Moreover, we prove
subsequently that under certain conditions on the similarity boundary, the SOSC
(restricted to $K$) is indeed satisfied. 
In order to highlight the usefulness of the similarity 
boundary in analyzing the properties of self-similar sets, we conclude with a 
few examples, chief among them the Hilbert space example found in [15].
Also, as we show later, in the case of self-similar 
tiles, the similarity boundary that we define, under certain conditions, 
does indeed coincide with the topological boundary.

\heading
\S2. Definitions and Notation
\endheading
Let $(X, d)$ be a complete metric space and
let $f_i: X \longrightarrow X$,  $i = 1, \cdots, N$ be
contracting similitudes with contraction factors $0 < r_i < 1$. Let $K \subset
X$ denote the
the unique non-empty, compact subset of $X$ such that 
$K = \bigcup_{i=1}^N f_i(K)$. $K$  is called the self-similar set or the 
invariant set of the IFS. We denote by $D$ the Hausdorff metric on the space
of non-empty, compact subsets of $X$. The Hausdorff metric is defined
the following way: $$D(F, K) = \inf \{ \epsilon : U(F, \e) \supset K 
\text{ and } U(K, \e) \supset F \}$$  where $U(A, \e)$  
denotes the $\e$-neighbourhood of $A$.

We define $\Cal I_n = \{ I = (i_1, \dots, i_n) : |I| = n \}$ 
for $n \geq 0$ where the length of $I$ is denoted by $|I|$ and $i_k \in
\{1, \dots, N\}$. For $n = 0$ the empty string $\emptyset$ is assumed to be 
an element of the set. 
We let $\Cal I = \bigcup_{n=0}^\infty \Cal I_n$.

We denote the concatenation of two strings $I = (i_1, \cdots, i_n)$ and
$J = (j_1, \cdots, j_m)$ the following way:
$IJ = (i_1, \cdots, i_n,j_1, \cdots, j_m)$.

$I$ and $J$ are termed incomparable if there exists no $L$ such that $I = JL$
 or $J = IL$. We denote incomparable strings the following way: 
$I \not\gtrless J$. The following is a list of short-hand notation that will
be used in the paper:
$$\align
f_I(K) &= f_{i_1}\circ \cdots \circ f_{i_n}(K)\\
K_I &= f_I(K)
\endalign$$
$$
r_I = \cases r_{i_1} \cdots r_{i_n} &\text{for }n>0\\
    1 &\text{for }n=0 \endcases
$$
$$
r_I^* = \cases r_{i_1} \cdots r_{i_{n-1}} &\text{for } n > 1\\
           1 &\text{for } n = 1\\
           \infty &\text{for } n = 0 \endcases
$$
$$
\text{diam}(K_I) = \cases r_I \t{diam}(K) &\text{for } I \neq \emptyset\\
              \t{diam}(K) &\text{for }I = \emptyset \endcases
$$
We denote by $r_{max}$ and $r_{min}$ the maximum and the minimum of the 
contraction factors.

The similarity dimension $\a$ is the unique real number satisfying the
equation $\sum_{i=1}^N {r_i}^\a = 1$.

The {\it sequence space} (or the {\it code space}) $\Omega$ is given by: 
$\Omega = \prod_{1}^\infty \{1, \cdots, N\}$.
Let $C_{i_1 \cdots i_n} = (i_1, \dots, i_n)\times
\prod_{n+1}^\infty \{1, \cdots, N\}$
denote the cylinder set in $\Omega$ whose first $n$ coordinates are $(i_1,
\dots, i_n)$. We define a metric $\rho$ on $\O$ the following way:
$$
\rho((i_n), (j_n)) = \cases 1 &\text{ if } i_1 \neq j_1\\
 r_{i_1} \cdots r_{i_m} &\text{ if } i_k = j_k\,\,\, 1 \leq k \leq m 
                           \text{ and } i_{m+1} \neq j_{m+1}\\
 0 &\text{ if } i_k = j_k \,\,\, \forall k \endcases
$$

With respect to this $\rho$ on $\O$ the cylinder sets are open (and closed). 
Since the cylinder sets generate the Borel sets of $\O$, we define a measure 
on $\O$ by defining it on the cylinder sets and extending it to the whole 
of $\O$. On the cylinder sets, we require that the measure $\nu$ have the 
property 
$\nu(C_{i_1 \cdots i_n}) = {(r_{i_1} \cdots r_{i_n})}^\a$. 
Again, for this metric $\rho$, $\nu$ is the Hausdorff $\a$-measure on $\O$ and
has the property that 
$\Cal H^\a(\O) = \nu(\O) = 1$ and the Hausdorff dimension of $\O$ is 
$\text{dim}_H \O = \a$ [3].
Corresponding to the similitudes $f_i$ on $K$, we define the similitudes
$\overline{f_i}$ on $\O$  by $\overline{f_i}((i_1, i_2, \cdots)) = (i, i_1, i_2, \cdots)$. 
They have the same corresponding contraction factors $r_i$.
Let the function $g:\O \longrightarrow K$ map the following way:
$g((i_1, i_2, \cdots)) =  
\lim_{n\to \infty}f_{i_1}\circ f_{i_2} \circ \cdots \circ f_{i_n}(K)$.
Then $g$ is a well-defined, continuous, Lipschitz map with Lipschitz constant
$\text{diam}(K)$. We define the image measure $\mu$ on $K$ through the map $g$.
Let $\mu(A) = \nu(g^{-1}(A)) \text{ for } A \subset K$. The measure $\mu$ has
been widely studied and seems to have its roots in the folklore of this 
subject. It plays a significant role in [1] and [14]. Since $g(C_i) = K_i$ but
$g^{-1}(K_i) \supset C_i$, we find that $\mu(K_i) \geq {r_i}^\a$. In fact
$\mu(K_i) = {r_i}^\a$ is equivalent to the condition $\mu(K_i \cap K_j) = 0$
for all $i \neq j$. It is also known that if $\mu(K_i \cap K_j) = 0$
for all $i \neq j$, then $\mu(K_I \cap K_J) = 0$ for all $I, J$ incomparable.
Similarly, if $\mu(K_i) = {r_i}^\a$ for $i \in \{1, \cdots, N\}$, then 
$\mu(K_I) = {r_I}^\a$ for all $I \in \Cal I$.

We now define the similarity boundary of the set $K$ denoted by $\partial_sK$.
$$B = \partial_sK = \bigcup_{j \neq k}{f_j}^{-1}(f_j(K) \cap f_k(K))$$

For example in the case of the Koch curve (Example 4.2)
 in ${\Bbb R}^2$, the boundary  
$\partial_sK$ is just 
the two endpoints. This is not the same as the topological boundary of the 
Koch curve which the curve itself.  On the other hand, the von Koch curve
is homeomorphic to an interval and the endpoints are its boundary as a
manifold.  So, even in this case there is some justification in the
analogy of $\partial_sK$ being its boundary.

Since $\partial_sK$ is closed, the set $U$ defined
as $U = K \setminus B$ is open in $K$. As we will see later $U$ can also be 
defined another way, independent of the definition of $\partial_sK$. We say 
that the similarity boundary $B$ is inverse invariant if 
${f_i}^{-1}(B)\cap K \subseteq B$ for all $i$. 

Lastly, we define the {\it strong open set condition in K } 
$\overline{\t{SOSC}}_K$. If there exists a set $U$ open in 
$K$ such that $f_i(U) \sbsq U$ for all $i$ and $f_i(U) \cap K_j = \es$ for all
$i \neq j$, then we say that the IFS satisfies the 
$\overline{\t{SOSC}}_K$. We would like to point out here that while we are
weakening the SOSC by requiring only that $U$ be open in $K$, we are also
requiring that the images of the open set $U$ be disjoint from the images of
the set $K$ itself. But as long as we restrict ourselves to $K$, this is only
a slightly stronger version of the SOSC.  Moreover, when the similitudes are
homeomorphisms of the ambient space as would be the case if $X$ is ${\Bbb R}^n$, 
then SOSC implies $\overline{\t{SOSC}}_K$.

\heading
\S3. Main Results 
\endheading

In this section, we use $B$ to refer to the similarity boundary and
$U$ to refer to $K \sm B$.

\proclaim{Lemma 3.1}${f_i}^{-1}(B) \cap K \subseteq B$ if and only if 
$f_i(U) \subseteq U$.
\endproclaim

\demo{Proof} ($\Rightarrow$) Let ${f_i}^{-1}(B) \cap K \subseteq B$.
Then, 
$$K \sm {f_i}^{-1}(B) \supseteq K \sm B = U$$
That is,
$$f_i(K) \sm B \supseteq f_i(U)\tag1$$
Since $ U = K \sm B\;$, $U \supset f_i(K)\sm B$. So $U \supset f_i(U)$ by (1).
\medskip

($\Leftarrow$) Conversely, let $f_i(U) \subseteq U$.
That is, $f_i(U) \subseteq K \sm B$.

So $$K \sm f_i(U) \supseteq B$$
$${f_i}^{-1}[K \sm f_i(U)] \cap K \supseteq {f_i}^{-1}(B) \cap K$$
$$[{f_i}^{-1}(K) \sm U] \cap K \supseteq {f_i}^{-1}(B) \cap K$$
$${f_i}^{-1}(K) \cap K \sm U \supseteq {f_i}^{-1}(B) \cap K$$
Since $B = K \sm U\;$, $$B \supset {f_i}^{-1}(K) \cap K \sm U $$
This implies $B \supset {f_i}^{-1}(B) \cap K$.
\qed
\enddemo

\proclaim{Lemma 3.2} $U = K \sm B = {\bigcup}_{i=1}^N {f_i}^{-1}(U_i) \cap K
\;$ where, $U_i = K \sm {\bigcup}_{i \neq j}f_j(K)$.
\endproclaim

\demo{Proof} 
$$\align
{f_i}^{-1}(U_i) \cap K &= {f_i}^{-1}[K \sm \bigcup_{j \neq i}f_j(K)] \cap K \\
&= [{f_i}^{-1}(K) \sm \bigcup_{j \neq i}{f_i}^{-1}(f_j(K))] \cap K \\
&= {f_i}^{-1}(K) \cap K \sm \bigcup_{j \neq i}{f_i}^{-1}(f_j(K)) \cap K  \\
&= K \sm \bigcup_{j \neq i}{f_i}^{-1}(f_j(K)) \cap K  \\
\endalign$$

Now taking the union over all the $i = 1, \cdots, N$,
$$\align
\bigcup_{i=1}^N {f_i}^{-1}(U_i) \cap K &=  K \sm \bigcup_{i=1}^N \bigcup_{j 
\neq i}
{f_i}^{-1}(f_j(K)) \cap K  \\
&= K \sm B
\endalign$$
This proves the Lemma. \qed
\enddemo

\proclaim{Lemma 3.3} $f_i(U) \cap f_j(K) = \emptyset$ for all $i \neq j$
\endproclaim

\demo{Proof}
$$\align
f_i(U) \cap f_j(K) &= f_i(K \sm B) \cap f_j(K) \\
&= [f_i(K) \sm f_i(B)] \cap [f_j(K)] \\
&= [f_i(K) \cap f_j(K)] \sm [f_i(B)]
\endalign$$

By the definition of $B$, 
$$\align
f_i(B) &= f_i \left(\bigcup_{j \neq k}{f_j}^{-1}(f_j(K) \cap f_k(K))\right) \\
&\supseteq f_i(K) \cap f_j(K)
\endalign$$

So, $f_i(U) \cap f_j(K) = \emptyset \;\; \forall i \neq j$. \qed
\enddemo

Thus, with the above result, if the set $U$ is non-empty, it satisfies the 
$\overline{\t{SOSC}}_K$. 

\proclaim{Lemma 3.4}Let $\mu(K_i) = {r_i}^{\a} \;\;\; \forall i=1, \cdots, N$. 
Then for any $\mu$-measurable 
set $A \subseteq K$, $\mu({f_I}(A)) = {r_I}^{\a}{\mu(A)}$ for any 
$I \in \Cal I$.

\endproclaim
 
\demo{Proof}
Let $A \subseteq K$ be a $\mu$-measurable set. For any $I \in \Cal I $,
consider the following commutative diagram:

$$
\CD
g^{-1}(A)     @>g>>    A  \\
@V\overline{f_I}VV      @VVf_IV  \\
\overline{f_I}(g^{-1}(A))     @>>g>   f_I(A)
\endCD
$$   
$$\mu(f_I(A)) = \nu(g^{-1}(f_I(A))) \geq \nu(\overline{f_I}(g^{-1}(A))) $$
where the $\overline{f_I}$ are the corresponding similitudes defined on 
$\Omega$.

Since on $\Omega$, $\nu$ is the Hausdorff measure:

$$\nu(\overline{f_I}(g^{-1}(A))) = {r_I}^{\a}\nu(g^{-1}(A)) = {r_I}^\a\mu(A)$$

So,$$\mu(f_I(A)) \geq {r_I}^{\a}\mu(A)\tag1$$

Since $A \subseteq K$,
 $$K_I = (K_I \sm f_I(A)) \cup f_I(A)$$
$$K_I \sm f_I(A) = f_I(K \sm A)$$

So,$$K_I = f_I(K \sm A) \cup f_I(A)$$
$$\mu(K_I) = \mu(f_I(K \sm A)) + \mu(f_I(A))$$
since the union is disjoint.

That is,$${r_I}^\a = \mu(f_I(K \sm A)) + \mu(f_I(A))$$

By (1):$${r_I}^\a \geq {r_I}^\a\mu(K \sm A) + {r_I}^\a\mu(A) = {r_I}^\a
\mu((K \sm A) \cup A) = {r_I}^\a$$

Hence,$$\mu(f_I(A)) = {r_I}^\a\mu(A)$$ 
which proves the Lemma.\qed
\enddemo

We are now ready to prove the main result.

\proclaim{Theorem 3.5}Let $({f_i}^{-1}(B) \cap K) \subseteq B \; \forall i$. That is,
let the similarity boundary be inverse invariant. Then the following are 
equivalent:
\roster
\item $\overline{\t{SOSC}}_K$ is satisfied by $\{f_i\}_{i=1}^N$
\item ${\t{\rm dim}}_H K = \a$
\item $U \neq \emptyset$
\item ${\t{\rm\i}} B = \emptyset$   
\item $\mu(K_i) = {r_i}^\a \;\;\; i=1, \cdots, N$
\item $\mu(K_i \cap K_j) = 0$
\item $\mu(B) = 0$
\endroster

\endproclaim

\demo{Proof}

(1) $\implies$ (2):

The following proof is a slight modification of the proof to be found in [15]. 
The hypotheses have been changed by doing away with the condition that the 
similitudes be onto and also that the SOSC be satisfied. Here, we assume only
that the $\overline{\t{SOSC}}_K$ is satisfied and the collection of similitudes
is not assumed to be bijective.
Let $U$ be open in $K$ such that $U$ satisfies all the rquirements of the 
$\overline{\t{SOSC}}_K$. Also let $\beta < \a$ be the Hausdorff dimension of 
$K$.
Let $I \in \Cal I$ be such that $K_I \sbsq U$. So for
every integer $n \geq 0$, the sets $\{K_{JI}\}_{|J| = n}$ are pairwise 
disjoint. Hence the self-similar set corresponding to the IFS 
$\{f_{JI}\}_{|J| = n}$ has disjoint images under the similitudes which generate
it. If the similarity dimension of this self-similar set is assumed to be
${\a}_n$, then its Hausdorff dimension will also be the same. Thus, 
 $${\a}_n \leq \t{dim}_H K \leq \a$$
Since $\a_n$ is the similarity dimension
$$1 = {\sum}_{|J|=n} {r_{JI}}^{\a_n} 
= {r_I}^{\a_n} {\sum}_{|J|=n}{r_{J}}^{\a_n}$$
So we get
$$\align
 {r_I}^{- {\a_n}} &= {\sum}_{|J|=n}{r_{J}}^{\a_n} \\
&\geq {\sum}_{|J|=n}{r_{J}}^\beta \\
&= {\sum}_{|J|=n}{r_{J}}^\a {r_{J}}^{\beta - \a} \\
&\geq {\sum}_{|J|=n}{r_{J}}^\a {r_{max}}^{n(\beta - \a)} \\
&=  {r_{max}}^{n(\beta - \a)}
\endalign $$

But on the one hand, while ${r_{max}}^{n(\beta - \a)}$ goes to infinity as n 
goes to infinity, on the other ${r_I}^{- \a_n}$ is bounded by ${r_I}^{-\beta}$.
This contradiction implies that the Hausdorff dimension of $K$ is $\a$.

\medskip

(2) $\implies$ (3):

If $U = \es$, then $K = B$ and so, $\i B \neq \es$. Then for $x \in B$, 
there exists $x \in G \sbs B\;\;\;$ open in $K$.
$$B = \bigcup_{i \neq j} {f_i}^{-1}(f_j(K)) \cap K$$

By the Baire Category Theorem, there exists $O \subseteq G$ open, such that,
$$O \subseteq {f_i}^{-1}(f_j(K)) \cap K$$

Since $ g^{-1}(O) \sbsq \O$, is open in $\O$, and the cylinder sets form a 
basis for the Borel sets in $\O$, 
there exists an $I \in \Cal I$ such that for some cylinder set $C_I$,
$$C_I \sbsq g^{-1}(O) \sbsq g^{-1}({f_i}^{-1}(f_j(K)) \cap K)$$

$$\implies g(C_I) = K_I \sbsq {f_i}^{-1}(f_j(K)) \cap K$$
 
That is,
$$f_i(K_I) \sbsq f_j(K) \cap f_i(K)$$
 
So $K_{iI} \sbsq K_j\;\;$ where clearly $iI \not\gtrless j$.
This means that the Hausdorff dimension $\t{dim}_H K < \a$ which is a 
contradiction.

\medskip

(3) $\implies$ (4):

Let $U \neq \es$ and $\i B \neq \es$. We have proved, in (2) $\implies$ (3)
that $\i B \neq \es$ implies that $\t{dim}_H K < \a$.
On the other hand, $({f_i}^{-1}(B) \cap K) \sbsq B \implies f_i(U) \sbsq U$ where, by
assumption, $U \neq \es$. So $U$ satisfies the $\overline{\t{SOSC}}_K$. 
This implies, from (1) $\implies$ (2)
that  $\t{dim}_H K = \a$, which is a contradiction. So $\i B = \es$.

\medskip

(4) $\implies$ (1):

Let $\i B = \es$. If $U = \es$, then $K = B$ and that is contrary to the 
assumption. So $U \neq \es$. So $U$ satisfies $\overline{\t{SOSC}}_K$.

\medskip

(1) $\implies$ (5):

\proclaim{Claim} $ g(\overline{f_i}(A)) = f_i(g(A))$ for all $A \sbsq \O$.
\endproclaim

Let $(i_1, i_2, \cdots ) = x \in A$. Then,
$$\align
\overline{f_i}(x) &= (i, i_1, i_2, \cdots ) \\
g(\overline{f_i}(x) &= \lim_{n \to \infty} 
f_i\circ f_{i_1} \circ \cdots \circ f_{i_n}(K) \\ 
&= f_i(\lim_{n \to \infty}f_{i_1} \circ \cdots \circ f_{i_n}(K)) \\
&= f_i(g(x))
\endalign$$
Since $x \in A$ is arbitarary, $ g\circ\overline{f_i} = f_i\circ g$.
\medskip
Now assume that $U \neq \es$, that $f_i(U) \sbsq U$ for all $i$ and 
that $f_i(U) \cap f_j(U) = \es$ for all $i \neq j$.
Recall that $U = K \setminus B$ with $B = \partial_sK$.  Thus,
$U$ is open in $K$.
Then, $\bigcup_{i=1}^N f_i(U) \sbsq U \sbsq K$.

Consider, 
$$\align
\s^{-1}(g^{-1}(U)) &= \bigcup_{i=1}^N \overline{f_i}(g^{-1}(U)) \\
g(\s^{-1}(g^{-1}(U))) &= g(\bigcup_{i=1}^N \overline{f_i}(g^{-1}(U))) \\
&= \bigcup_{i=1}^N f_i(U) \t{  by Claim} \\
&\sbsq U
\endalign$$
 
This implies $\s^{-1}(g^{-1}(U)) \sbsq g^{-1}(U)$. So $g^{-1}(U)$ is invariant
under the inverse of the Bernoulli shift map.   

Therefore, by the Ergodic Theorem, $\nu(g^{-1}(U)) = 1$. That is, $\mu(U) = 1$.
But, since $\nu$ is the Hausdorff measure on $\O$ and $\overline{f_i}$ is a
similitude,
$$ \nu(\overline{f_i}(g^{-1}(U))) = {r_i}^\a \nu(g^{-1}(U)) = {r_i}^\a\tag2$$

By the Claim,
$$g(\overline{f_i}(g^{-1}(U))) = f_i(U)$$
So,
$$\overline{f_i}(g^{-1}(U)) \sbsq g^{-1}(f_i(U))$$ 
 
From (2),
$$\mu(f_i(U)) = \nu (g^{-1}(f_i(U))) 
\geq \nu(\overline{f_i}(g^{-1}(U))) = {r_i}^\a$$

Since $\{f_i(U)\}_{i=1}^N$ is a disjoint collection,

$$1 = \mu(U) \geq \mu(\bigcup_{i=1}^N f_i(U)) = \sum_{i=1}^N \mu(f_i(U))
\geq \sum_{i=1}^N  {r_i}^\a = 1$$

This forces $\mu(f_i(U)) =  {r_i}^\a$ and 
$\mu(U) = \mu(\bigcup_{i=1}^N f_i(U))$.

Now consider,
$$\mu(f_1(K)) + \mu(f_2(K)) + \cdots + \mu(f_N(K)) \geq 1$$
 
But $f_1(K) \cap f_j(U) = \es$ for all $j \neq 1$. So we have,
$$\align
1 = \mu(K) &\geq \mu( f_1(K) \bigcup (\bigcup_{i=2}^N f_i(U))) \\
&= \mu(f_1(K)) + \sum_{i=2}^N \mu(f_i(U)) \geq 1
\endalign$$
This forces $\mu(f_1(K)) = {r_1}^\a$. The same argument can be repeatedly
used to show that $\mu(f_i(K)) = {r_i}^\a \;\;\; \forall i$.

\medskip

(5) $\implies$ (6):

Let $\mu(K_i) = {r_i}^\a \;\;\; i=1, \cdots, N$. 
$$\nu(g^{-1}(K_i) \sm C_i) = \nu(g^{-1}(K_i)) - \nu(C_i) = 0$$

$g^{-1}(K_i)$ may be written in the following way,
$$g^{-1}(K_i) = (g^{-1}(K_i) \sm C_i) \cup C_i$$
It follows that for any $i \neq j$,
$$\nu(g^{-1}(K_i) \cap C_j) = \nu((g^{-1}(K_i) \sm C_i) \cap C_j) \cup 
\nu(C_i \cap C_j) = 0$$
So,
$$\align
\mu(K_i \cap K_j) &= \nu(g^{-1}(K_i) \cap g^{-1}(K_j)) \\ 
&= \nu((g^{-1}(K_i) \sm C_i)\cap g^{-1}(K_j)) + \nu(C_i \cap  g^{-1}(K_j)) = 0
\endalign$$

\medskip

(6)  $\implies$ (5):

Let $\mu(K_i \cap K_j) = 0$. That is, $\nu(g^{-1}(K_i) \cap g^{-1}(K_j)) = 0$.
Since $\mu(K_i) = \nu(g^{-1}(K_i)) \geq {r_i}^\a$, we need to show that
$\mu(K_i) \leq {r_i}^\a$.

Consider,
$$\align
g^{-1}(K_i) &= C_i \bigcup \left[ \bigcup_{j \neq i} 
g^{-1}(K_i) \cap C_j\right] \\
&\sbsq C_i \bigcup \left[\bigcup_{j \neq i}g^{-1}(K_i) \cap g^{-1}(K_j)\right]
\endalign$$

By hypothesis $\nu(\bigcup_{j \neq i}g^{-1}(K_i) \cap g^{-1}(K_j)) = 0$.
So, $\mu(K_i) = \nu(g^{-1}(K_i)) \leq {r_i}^\a$.

\medskip

(6) $\implies$ (7):

Let $\mu(K_i \cap K_j) = 0$.
$$B = \bigcap_{i \neq j} {f_i}^{-1}(K_i \cap K_j)$$
Let ${f_i}^{-1}(K_i \cap K_j) \cap K = A$. Then $K_i \cap K_j = f_i(A)$.

In order to prove that $\mu(B) = 0$, we need only prove that $\mu(A) = 0$.
Now $\mu(f_i(A)) = 0$ by hypothesis. But $0 = \mu(f_i(A)) = {r_i}^\a \mu(A)$ 
by Lemma 3.4.
So $\mu(A) = 0$ and hence  $\mu(B) = 0$.

\medskip

(3)  $\implies$ (7):

$U \neq \es$ and $f_i(U) \sbsq U$ implies, by the Ergodic Theorem, that
$\mu(U) = 1$. Therefore $\mu(B) = \mu(K \sm U) = 0$.

\medskip

(7)  $\implies$ (3):

$\mu(B) = 0$ implies that $\t{int} B = \es$. But if $U = \es$, then $K = B$
which is a contradiction. So $U \neq \es$.
 
\enddemo

\bigskip

{\bf Remark.} The proof of (1) implies (5), that is, the
$\overline{\t{SOSC}}_K$
implying that for all $i \neq j$
$\mu(K_i \cap K_j) = 0$ 
did not use the invariance of the
similarity boundary and as such stands by itself. This is also true of 
the proof of the 
equivalence of (5) and (6). It is possible to prove Theorem 3.5 by proving
fewer implications than above. We have included more
proofs than necessary for
the insights they provide.

\bigskip

We now prove that in the case of self-similar tiles in Euclidean space, when
the similarity boundary is inverse invariant, it is the same as the topological
boundary. Intuitively, this is what we would expect and also what would be a 
desirable property in the similarity boundary.

\bigskip

\proclaim{Theorem 3.6} Let $K$ be a self-similar tile in n-dimensional 
Euclidean space ${\Bbb R}^n$ with similarity dimension $n$. If the 
similarity boundary of $K$ is inverse invariant, then
$$\partial_sK = \partial K$$ 
where $\partial K$ represents the topological boundary.
\endproclaim

\demo{Proof}
 
We first prove that $\partial_sK \sbsq \partial K$. Since $K$ is a tile in 
${\Bbb R}^n$, it has non-empty interior and since the similarity dimension of 
$K$ is $n$, the Hausdorff dimension is also $n$. Let $\t{int}K = O$.

\medskip

{\bf Claim.} $\t{int}(f_i(K)) \cap \t{int}(f_j(K)) = \es$ for all $i \neq j$.

\medskip

{\bf Proof of Claim.} Since the similitudes are open maps in Euclidean space,
$\t{int}(f_i(K)) = O_i$. So we need to show that $O_i \cap O_j = \es$ for all
$i \neq j$. Let us assume the contrary. Then, since $O_i \cap O_j$ is an open
set, $\lambda (O_i \cap O_j) > 0$, where $\lambda$ is the Lebesgue measure in
${\Bbb R}^n$. But, from [1]:
$$0 = \lambda (K_i \cap K_j) \geq \lambda (O_i \cap O_j) > 0$$
which is a contradiction that proves the Claim.

$O_i \cap O_j = \es$ implies that 
$$K_i \cap K_j = \partial K_i \cap \partial K_j\tag1$$
But, since $f_i$ is a homeomorphism on ${\Bbb R}^n$, 
${f_i}^{-1}(\partial K_i) = \partial K$. This means that 
$\partial K_i = f_i(\partial K)$. So from (1):
$$K_i \cap K_j = f_i(\partial K) \cap f_j(\partial K)\tag2$$
Now, in order to prove that 
$\partial_sK \sbsq \partial K$, we need to prove that for every $i \neq j$,
${f_i}^{-1}(K_i \cap K_j) \sbsq \partial K$. That is, we need to prove that
for every $i \neq j$, $K_i \cap K_j \sbsq f_i(\partial K)$ which follows from
(2). So $\partial_sK \sbsq \partial K$.

We now prove the containment in the other direction. That is, 
$\partial K \sbsq \partial_sK$. $\partial_sK \sbsq \partial K$ is a closed set.
Let $x \in \partial K \sm \partial_sK$.

\medskip

{\bf Claim.} Let $x \in \partial K \sm \partial_sK$. For all $i$, 
$f_i(x) \in \partial K$.

\medskip

{\bf Proof of Claim.} If $f_i(x) \in \t{int}K$, then, since 
$f_i(x) \in \partial K_i$, $f_i(x) \in f_j(K)$ for some $j \neq i$. 
By the definition of the similarity boundary this would mean that 
$x \in \partial_sK$ and that is a contradiction. So for all
$i$, $f_i(x) \in \partial K$. This proves the Claim.

But on the other hand, if $f_i(x) \in {\partial}_sK$, then 
$x \in {f_i}^{-1}(\partial_sK) \sbsq \partial_sK$ by the assumption on the 
similarity boundary. But this contradicts the choice of $x$. Since 
$x \in \partial K \sm \partial_sK$ was arbitarary, this proves that 
$f_i(\partial K \sm \partial_sK) \sbsq \partial K \sm \partial_sK$ for
all $i$. So $f_i(\overline{\partial K \sm \partial_sK}) \sbsq
\overline {\partial K \sm \partial_sK}$ for all $i$. Let 
$A = \overline {\partial K \sm \partial_sK}$. So 
$\overline {\partial K \sm \partial_sK}$ is forward invariant. That is,
$\bigcup_{i=1}^N f_i(A) \sbsq A$. This implies that $K \sbsq A$ which is a 
contradiction by the Theorem 3.5 above. So $\partial K = \partial_sK$.

\enddemo

\heading
\S4. CONCLUSION
\endheading

It is worthwhile to point out that in Theorem 3.6 above, 
$\partial_sK \sbsq \partial K$ is true without the similarity boundary being
necessarily inverse invariant. However, the inverse containment is not always
true if the similarity boundary is not inverse invariant. We illustrate this
point with the following example.

\example{Example 4.1} In ${\Bbb R}^2$, we define the following four 
similitudes: $f_1$
is contraction towards $(0, 1)$ by a factor $1/2$, $f_2$ is contraction towards
$(1, 1)$ by $1/2$ , $f_3$ is contraction towards $(0, 0)$ by $1/2$ and $f_4$ is
contraction towards $(1, 0)$ by $1/2$. Then the self-similar set we get is 
the unit square in ${\Bbb R}^2$. It can be shown that its similarity boundary
is the outer edges of the square and that it is inverse invariant. It also
coincides with the topological boundary. Now, we modify the similitudes 
slightly the following way: we add a rotation by $-90$ degrees to $f_2$, a 
rotation
by $90$ degrees to $f_3$ and a rotation by $180$ degrees to $f_4$. 
We still get the
same invariant set. But now, the similarity boundary is just the line segments
connecting $(1, 1)$ to $(1, 0)$ and $(0, 0)$ to $(1, 0)$. Also, the similarity
boundary is no longer inverse invariant and it no longer coincides with the 
topological boundary of the unit square.

\endexample
 
As a further point of interest, it has
been proved that the Hausdorff dimension of $\partial_sK$ is in fact strictly
less than the similarity dimension $\a$ of the set $K$ [9]. We would like to
mention that the complement of the similarity boundary in $K$ is not 
necessarily the
maximal open set (open in $K$) satisfying the requirements of the 
$\overline{\t{SOSC}}_K$. As an illustration, consider the following example. 

\example{Example 4.2}
The von Koch curve is the invariant set of the IFS consisting of the
following contraction
similitudes $f_1, f_2, f_3$ and $f_4$ defined the following way.  Let $f_1$ be
contraction towards the origin $(0, 0)$ by a factor 1/3.  Let $f_2$ be contraction
by a factor $1/3$ towards the origin, rotation by $60$ degrees followed
by translation to the right by $1/3$ units so that $(0, 0)$ goes to $(1/3, 0)$. 
Let $f_3$ be contraction
by a factor $1/3$ towards the origin, rotation by $-60$ degrees followed
by a translation taking $(0, 0)$ to $(\frac12, \frac{\sqrt{3}}{6})$.  Let
$f_4$ be
contraction towards the point $(1, 0)$ by a factor 1/3. Let $U$
be the set formed by deleting the origin from the curve. Then $U$ is open in
$K$.  It can be easily verified that it contains the complement of the 
similarity boundary.  Furthermore, it satisfies the requirements of the
$\overline{\t{SOSC}}_K$.

\endexample

Lastly, we take another look at an example in [15] of a self-similar set $K$ 
in Hilbert space which has the property that it satisfies the SOSC but has
zero Hausdorff $\a$ measure. The example that follows is basically the one
in [15]. However, in our analysis, there is no need for the similitudes to be 
onto and so the definitions are somewhat simpler.

\example{Example 4.3}
We consider the $l_1$ space endowed with the corresponding norm. 
$$l_1 = \{(x_n): x_n \in {\Bbb R}, \sum_{n=1}^{\infty} |x_n| < \infty\}$$
Let $\{e_n\}_{n=1}^{\infty}$ denote the unit vectors which form a basis for the
space. We define the similitudes $f_1, f_2$ and $f_3$ the following way for
$x = (x_1, x_2, x_3, \cdots) \in l_1$,
$$f_1(x_1, x_2, x_3, \cdots) = (0, x_1/2, 0, x_2/2, 0, x_3/2, 0, \cdots)$$
$$f_2(x_1, x_2, x_3, \cdots) = (0, 0,x_1/2, 0, x_2/2, 0, x_3/2, 0, \cdots)$$
$$f_3(x_1, x_2, x_3, \cdots) = ({(x_1 + 1)}/2, {x_2}/2, {x_3}/2, \cdots)$$
$f_1, f_2$ and $f_3$ are contracting similitudes on $l_1$ with contraction
factors $r_1 = r_2 = r_3 = 1/2$. $f_1$
and $f_2$ have as their fixed point the origin and since $f_3$ can be written
as $f_3(x) = e_1 + 1/2(x-e_1)$, the fixed point of $f_3$ is $e_1$. It is clear 
from the definition of $f_1$ and $f_2$ that $K_1 \cap K_2 = (0,0,0, \cdots)$. 
Also, if $f_3(x_1, x_2, x_3, \cdots) = f_1(y_1, y_2, y_3, \cdots)$, then
$x_1 = -1$. Now let $P$ be the set of points in $l_1$ with non-negative 
coordinates. Then, $f_i(P) \sbsq P$ for $i = 1,2,3$. Since $P$ is a closed set,
this implies that $K \sbsq P$ and so, $K$ cannot contain any points with
negative coordinates. This in turn implies that $K_3 \cap K_1 = \es$. The same
argument also shows that $K_3 \cap K_2 = \es$. So $f_i(K) \cap f_j(K)$ is 
either empty or the origin. Since the origin is the fixed point of both $f_1$
and $f_2$, ${f_1}^{-1}(0) = 0 = {f_2}^{-1}(0)$. ${f_3}^{-1}(0)$ is a point 
whose first coordinate is $-1$ and therefore lies outside $K$. All this proves
that the similarity boundary of this set $K$ is nothing but the origin and 
also that it is inverse invariant. Since $\mu(B) = 0$, all the equivalences
from Theorem 3.5 follow. So if we restrict ourselves to $K$, then instead of
having to explicitly construct an open set to satisfy the SOSC, Theorem 3.5 
gives us the existence of such a set. In particular, ${\t{\rm dim}}_H K = \a$
where $(1/2)^\a + (1/2)^\a + (1/2)^\a = 1$.
 Once again, restricting ourselves to $K$,
from [15], we get that $K$ has zero Hausdorff $\a$ measure. Hence we get the
required example. 

\endexample

In the above example, the contraction factors could have been any 
$\{c_1, c_2, \cdots, c_N\}$ with $0 < c_i < 1$ and the analysis would 
have been
the same. It is somewhat easier to reach the conclusions about $K$ than 
the analysis found in [15].

\enddocument